\documentclass[letterpaper, 10 pt, conference]{ieeeconf}  % Comment this line out if you need a4paper

\usepackage{amsmath}
\usepackage{amssymb}
\usepackage{graphicx}

\IEEEoverridecommandlockouts                              % This command is only needed if
\title{\LARGE \bf
Minimizing Risk of Load Shedding and Renewable Energy Curtailment in a Microgrid with Energy Storage
}

\author{Ashkan Zeinalzadeh$^{1}$ and Vijay Gupta$^{2}$% <-this % stops a space
\thanks{*This work was supported by the Center for Sustainable Energy at Notre Dame}% <-this % stops a space
\thanks{$^{1}$Ashkan Zeinalzadeh is with the Department of Electrical Engineering, University of Notre Dame, IN, USA
        {\tt\small azeinalz@nd.edu}}%
        \thanks{$^{2}$Vijay Gupta is with Faculty of Electrical Engineering, University of Notre Dame, IN, USA
        {\tt\small vgupta2@nd.edu}}%
}

\begin{document}

\maketitle
\thispagestyle{empty}
\pagestyle{empty}

%%%%%%%%%%%%%%%%%%%%%%%%%%%%%%%%%%%%%%%%%%%%%%%%%%%%%%%%%%%%%%%%%%%%%%%%%%%%%%%%
\begin{abstract}

We consider a microgrid with random load realization, stochastic renewable energy production, and an energy storage unit. The grid controller provides the total net load trajectory that the microgrid should present to the main grid and the microgrid must impose load shedding and renewable energy curtailment if necessary to meet that net load trajectory. The microgrid controller seeks to operate the local energy storage unit to minimize the risk of load shedding, and renewable energy curtailment over a finite time horizon. We formulate the problem of optimizing the operation of the storage unit as a finite stage dynamic programming problem. We prove that the multi-stage objective function of the energy storage is strictly convex in the state of charge of the battery at each stage. The uniqueness of the optimal decision is proven under some additional assumptions. The optimal strategy is then obtained. The effectiveness of the energy storage in decreasing load shedding and RE curtailment is illustrated in simulations.

\end{abstract}

%%%%%%%%%%%%%%%%%%%%%%%%%%%%%%%%%%%%%%%%%%%%%%%%%%%%%%%%%%%%%%%%%%%%%%%%%%%%%%%%
\section{INTRODUCTION}

As the penetration of Renewable Energies (RE) increases, electric power systems encounter new operating problems uncommon in conventional systems \cite{Ashkan1}. %\cite{Zeinalzadeh}-\cite{Ashkan}.
RE sources typically exhibit an intermittent pattern. This pattern presents significant challenges for utilities, e.g. for maintaining the transmission line capacity, regulating frequency, and balancing the generation and load \cite{Ashkan2}-\cite{Ashkan3}. %In addition, the level of intermittent renewable generation may exceed overall system capacity.
%Here, we study the problem of maintaining the capacity of line connecting the main grid to a microgrid that is equipped with loads, RE units and energy storage. long transmission line equipped with RE at the end of the grid is critically important for the grid controller. The transmission line capacity is proportional to the square of voltage and inversely proportional to the line impedance \cite{whitepaper}. The line impedance increases with transmission distance.
Indeed, an excess of renewable generation at the time of low loads and transmission constraints can lead to curtailment of renewable energies \cite{whitepaper}-\cite{Golden}, meaning that the grid controller does not allow the RE resources to inject power into the grid. Interestingly, this may be concomitant with load shedding that occurs if the grid cannot meet the load because of transmission line constraints or supply-demand imbalance.
%Transmission constraints have been the most common cause of wind energy curtailment in the United States \cite{NREL}. An efficient charging and discharging strategy for energy storage can enhance the transmission line capacity for a higher penetration of REs in remote grids, e.g. charging the batteries at times of excess energy, and discharging at times of peak load. %However, the uncertain and stochastic nature of the load and REs and
%The lack of communication among consumers creates an inability to coordinate the charging/discharging strategies of distributed energy storage units. Therefore, it is more viable and economical for consumers to invest in a large capacity energy storage unit that operates locally by community aggregator.

In this work, we consider a microgrid that is equipped with RE and energy storage (in the form of a battery). The energy storage unit is operated by the non-profit microgrid operator, who can
%which is a separate agent from the grid controller. The consumers in the remote grid community own the RE resources and share values that trend towards increasing the penetration of renewable energy. The consumers' renewable energy in excess of their loads are injected into the grid. The community aggregation unit does not purchase and sell the energy, but contracts with consumers to maintain the transmission line capacity. The grid controller can
direct load shedding and renewable energy curtailment if needed. The main grid provides the microgrid operator with a maximum net load, e.g. due to transmission line power constraints, that it is allowed to present to the main grid. The problem for the microgrid operator is to minimize both RE curtailment and load shedding while meeting the total net load constraint.
%to avoid violating the operating constraints and satisfy the power flow equations. The community aggregator decides on
It seeks to perform this by choosing the charging and discharging strategies for energy storage to minimize the risk of high amounts of load shedding and RE curtailment.

The chief contribution of this paper is the formulation and analysis of this problem as a multi-stage dynamic programming problem. Specifically, we use the Conditional
Value-at-Risk (CVaR) as a metric of risk for the load shedding and RE curtailment and prove the strict convexity of the objective function under some mild conditions. Thus the optimization problem can be efficiently solved.
%We formulate the community aggregator (energy storage unit) objective function as a multi-stage dynamic program and prove the strict convexity of the objective function under sufficient conditions.
A numerical analysis is presented for optimization of the energy storage for a stochastic model for the load and RE output. It is shown that the optimal strategy for the energy storage enables consumers to increase the penetration of renewable energies.
%creates circumstances under which, excess renewable energies from consumers can be injected into the grid and used locally to meet the load at peak-times. Thus, energy storage enables consumers to increase the penetration of renewable energies. %We optimize the energy storage operation over finite time stages, due to the limited lifetime of energy storage. The goal of this study is to develop a decentralized (without communication with grid controller) optimal strategy for the operation of the energy storage unit based on the stochastic models of load and REs.
%The REs are generated locally and in a long distance from the grid controller.
%We define an objective function for the energy storage to decrease the risk of high amounts of load-shedding and RE curtailment.
%The energy storage objective function is formulated to decrease the risk of high amounts of renewable energy curtailment and load shedding.

The optimal use of the energy storage has been investigated for different applications, e.g. frequency regulation \cite{Outdalov} and \cite{Mercier}, to mitigate the intermittency of the renewable sources \cite{Kanoria}, \cite{Baker}, \cite{Wang}, and \cite{Huan}, as well as peak shaving \cite{Even}, and spinning reserves \cite{Kottick}. Conventional optimal power flow without storage, decouples the optimization in different time periods \cite{Chandy}. In contrast, energy storage introduces correlation across time periods. The main challenges to optimize the energy storage operation is this correlation across time, constraints on the capacity and charge/discharge rate of the energy storage unit, and stochastic behaviour of load and renewable energies. The studies most related to our work are \cite{Anderson} and \cite{Lou}. Authors in \cite{Anderson} study the load shedding for the frequency regulation in a grid model with a deterministic load. Authors in \cite{Lou} formulate the optimal load shedding due to supply reduction. They maximize the operator profit for a grid model without energy storage and with a deterministic load. Unlike these works, we minimize the risk (as measured by CVaR) of high amounts of load shedding and RE curtailment in the presence of stochastic load and RE generation.

%We assume the one stage energy storage objective function is a general convex function in the charge/discharge rate.
%We prove that the multi-stage objective function is convex in the state of the charge of the battery at each stage. Therefore, the energy storage operation can be optimized using convex optimization. Finally, we assume that the load and renewable resource output are independent random processes with known distributions for a finite number of discrete time slots.

The rest of the paper is structured as follows. In Section~\ref{SectionII}, an energy storage objective function is modeled with the goal of decreasing the risk of high amounts of load shedding and RE curtailment in a microgrid. In Section~\ref{SectionIII}, we formulate the optimization problem of energy storage unit as a dynamic programming problem and prove the convexity of the objective function in action space at each stage. We derive a sufficient condition for the strict convexity of the objective function. In Section~\ref{SectionV} we show a numerical result obtained through simulation evaluations. Finally, concluding remarks are provided in Section~\ref{SectionVI}.\\

\section{MICROGRID MODEL}\label{SectionII}

\begin{figure}[!htb]
  \centering
  \includegraphics[width=8.5cm,height=7cm,keepaspectratio]{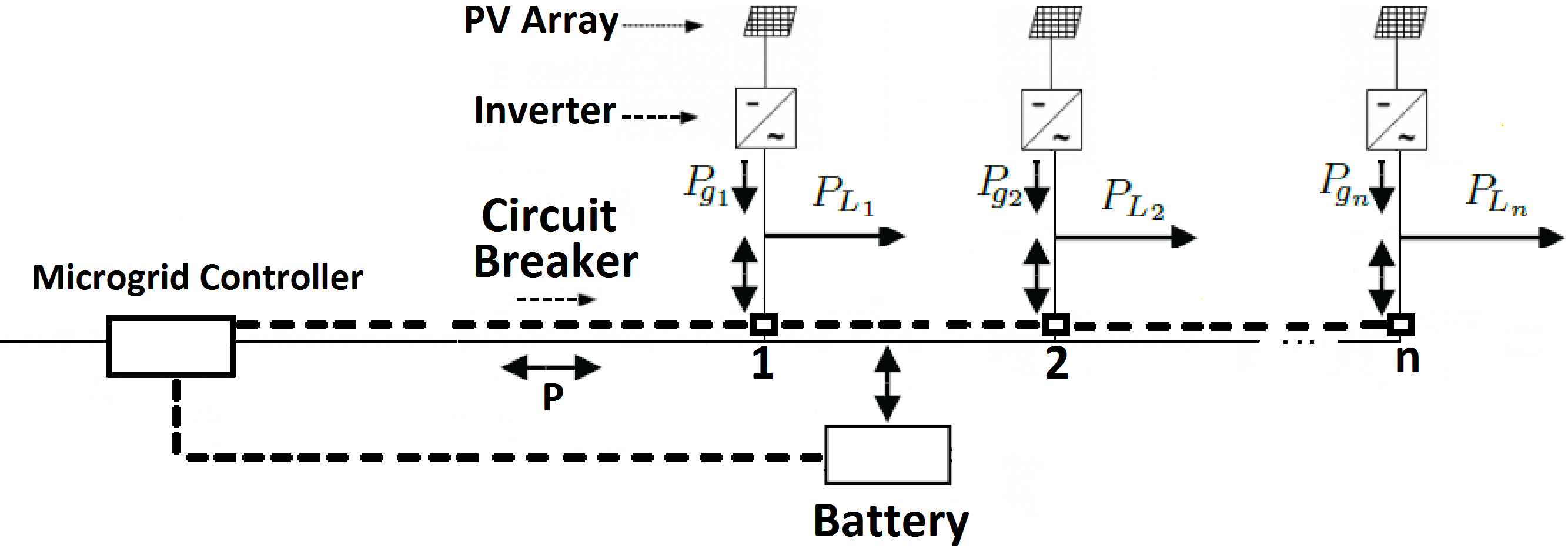}\\
  \caption{Microgrid model: The dashed lines are the control signals and the solid lines
are the power flow lines. The microgrid controller directs load shedding and RE
curtailment using circuit breakers to maintain transmission line constraint.}\label{fig:fig1}
\end{figure}

We consider a microgrid with local renewable energy resources (e.g. solar PV systems) and an energy storage unit as shown in Figure~\ref{fig:fig1}. Let $P_{L_k}$ and $P_{g_k}$ denote respectively the load and output of solar PV system of the $k$th consumer. %The dashed lines are the control signals and the solid lines are the power flow lines.
The microgrid controller allows the distributed renewable generators to provide necessary active power locally or to the main grid. In addition, it directs load shedding and RE curtailment to obey a constraint on the max total active power that the microgrid can draw from, or inject into, the transmission line that connects it to the main grid. To this end, it uses the energy storage unit to minimize the curtailment of renewable energies and to minimize the load shedding. %It is desirable to meet the load locally using RE resources and the energy storage unit, so as to relax the power flow on the long transmission line.
Intuitively, the energy storage unit can act as an energy reserve. %to prevent load shedding. We formulate the energy storage problem such that the renewable energy's integration into the remote grid becomes more reliable, i.e.
However optimizing the usage of energy storage to minimize the risk of high amounts of load shedding and RE curtailment is a difficult and open problem.

The following notations are used in this work:\\

\textbf{Notation:}
\begin{itemize}
\item $d_{t}$ and $r_{t}$ are the load and renewable resource realizations respectively at time $t$.
\item $\tilde{d}_{t}$ and $\tilde{r}_{t}$ are the load shedding and renewable resource curtailment respectively at time $t$.
\item $|.|$ denotes the absolute value.
\item $(.)^+=\max(.,0)$.% and $(.)^-=\min(.,0)$.
\end{itemize}

Electricity generated by renewable resources would first supply the local load. Any additional power needed to satisfy the load could be supplied from the grid or satisfied by discharging the battery. The energy in excess of the load, that is generated by renewable resources, can be used to charge the battery or be injected into the grid. The microgrid controller has the authority to curtail the excess RE injected into the grid and to direct the load to be shutdown. The total power to or from the microgrid must be less than a specific value that is specified by the main grid to satisfy, e.g, the transmission line constraints.
%that is assigned by the main grid controller. %(transmission line constraints).% specific value $P$.
%We consider a battery cost that is a function of charging and discharging power

\subsection{Load and Renewable Energy}

Let $\mathfrak{T}=\{1,2,...,T\}$ be the index set. Define $d=\{ d_{t},  t \in \mathfrak{T} \}$ and $r=\{ r_{t}, t \in \mathfrak{T} \}$ as random processes on the probability space $( \Omega,\mathfrak{F}, \mathbb{P})$. $d$ and $r$ represent the load and renewable resource power respectively. For a fixed $t$, and for $\omega \in \Omega$, $d_{t}(\omega)$ and $r_{t}(\omega)$ are non-negative random variables with known continuous probability density functions. For the given $\omega \in \Omega$, $d_{t}$, and $r_{t}$ are deterministic functions of $t$ that denote the load realization and renewable resource power realization respectively at time $t$. Let $P_{L_k}^t$ and $P_{g_k}^t$ be respectively the load and output of solar PV system of the $k$th consumer in Figure~\ref{fig:fig1} at time $t$. $d_t$ and $r_t$ are the aggregate consumers' load and renewable power realization respectively,
\begin{align}\label{eq:loadpv}
d_t=P_{L_1}^t+P_{L_2}^t+...+P_{L_n}^t, \\ \notag
r_t=P_{g_1}^t+P_{g_2}^t+...+P_{g_n}^t.
\end{align}
Below, the operating constraints for the energy storage are described.

\subsection{Energy Storage}

Energy storage can charge or discharge $b_{t}$ amount of power at time $t$, $b_{t}$ can be either positive (battery is charging) or negative (battery is discharging). Let $s_{t} \geq 0$ denote the energy level of the battery at the beginning of the time $t$ for all $t=1,...,T$. State of charge of the battery, $s_t$, evolves according to
\begin{align}\label{eq:one}
%&s_{1}=s_{0}\\ \notag
s_{t}=a \,\ s_{t-1}+ b_{t-1} \,\ \Delta t
\end{align}
where $s_{t-1}$ is the energy level of the battery at the beginning of the time $t-1$, that is reduced by the factor $0 < a < 1$ at time $t$. It is assumed that $b_t$ is constant in the time slot $[t,t+\Delta t)$. The energy level of the battery is bounded by a minimum and maximum capacity
\begin{align}\label{eq:two}
s^{\min} \leq s_{t} \leq s^{\max}.
\end{align}
Using (\ref{eq:one}), inequality (\ref{eq:two}) can be written as follows
\begin{align}\label{eq:constraint}
\frac{s^{\min}-a s_{t}}{\Delta t} \leq b_{t} \leq \frac{s^{\max}-a s_{t}}{\Delta t}.
\end{align}
Maintaining constraints (\ref{eq:one}) and (\ref{eq:constraint}) is crucial in optimization of energy storage operations. The method by which the microgrid controller decides on the value of renewable energy curtailment and load shedding to maintain the transmission line constraint is described below.

\subsection{Renewable Energy Curtailment and Load Shedding}

Renewable resources can be used to meet the load or charge the battery locally. The energy in excess of the load that is generated from the renewable resources can be injected into the grid. The microgrid controller can prohibit renewable resources from injecting energy into the grid, a process that is called resource curtailment. The resource curtailment, $\tilde{r}_{t}$, is bounded by
\begin{align}\label{eq:ineqrc}
0 \leq \tilde{r}_{t} \leq (r_{t}-b_{t}-d_{t})^+.
\end{align}
The load can be satisfied by using renewable resources or discharging the battery locally, and the extra energy needed to meet the load can be supplied from the grid. To prevent the violation of the transmission line constraint, excessive load can be avoided in a process known as load shedding. The load shedding, $\tilde{d}_{t}$, is bounded by
\begin{align}\label{eq:ineqls}
0 \leq \tilde{d}_{t} \leq (d_{t}+b_{t}-r_{t})^+.
\end{align}
The power flow on the transmission line can be written as
\begin{align}\label{eq:power}
p_{t}=d_{t}-\tilde{d}_{t}+b_{t}-r_{t}+\tilde{r}_{t}.
\end{align}
The grid controller provides $P^{\text{min}}$ and $P^{\text{max}}$ to the microgrid controller as the lower and upper bounds on the power flow $p_t$
\begin{align}\label{eq:ineqp}
P^{\text{min}} \leq p_{t} \leq P^{\text{max}}.
\end{align}
The main grid operator controls the output power of the generator, and provides the microgrid
operator with a maximum net load. The microgrid operator controls the load shedding, $\tilde{d}_{t}$, and resource curtailment, $\tilde{r}_{t}$. The main grid controller's objective is to control the power flow to avoid violation of the operating constraints (\ref{eq:ineqrc})-(\ref{eq:ineqp}). Let $n_{t}=d_{t}-r_{t}$ and $\tilde{n}_{t}=\tilde{d}_{t}-\tilde{r}_{t}$. Inequalities (\ref{eq:ineqrc}) and (\ref{eq:ineqls}) can be written as
\begin{align}\label{eq:ineq11}
-(-n_{t}-b_{t})^+  \leq \tilde{n}_{t} \leq (n_{t}+b_{t})^+.
\end{align}
Equation (\ref{eq:power}) can be written as
\begin{align}\label{eq:newpower}
&p_{t}=n_{t}-\tilde{n}_{t}+b_{t}.
\end{align}
The microgrid controller decides on the value of $\tilde{n}_{t}$, which is given as a function of $b_t$ as follows
\begin{align}\label{eq:shedding2}
\tilde{n}_{t}^*(b_t)=\left\{ \begin{array}{rcl}
   n_{t}+b_{t} - P^{\max}   & \mbox{if} & (n_{t}+b_{t}) > P^{\max}\\
   n_{t}+b_{t} - P^{\min}   & \mbox{if} & (n_{t}+b_{t}) < P^{\min} \\
0 & \mbox{otherwise} & .
\end{array}\right.
\end{align}
Let $\tilde{\mathfrak{n}}_{t}(b_t)= \mid \tilde{n}_{t}^*(b_t) \mid$. We use the abbreviation $\tilde{\mathfrak{n}}_{t}$ instead of $\tilde{\mathfrak{n}}_{t}(b_t)$. Let $F_{\tilde{\mathfrak{n}}_{t}}$ be the cumulative distribution function for $\tilde{\mathfrak{n}}_{t}$. Below, the objective function of energy storage is defined in order to decrease the risk of high amounts of load shedding and RE curtailment.

\subsection{Risk Minimization}

To define the objective function of the energy storage, we use the notion of Value-at-Risk (VaR) and Conditional Value-at-Risk (CVaR) \cite{Uryasev}. $VaR_{\alpha}(\tilde{\mathfrak{n}}_{t})$ determines the worst possible $\tilde{\mathfrak{n}}_{t}$ that may occur with confidence level $\alpha$. For a given $0 < \alpha < 1$, the amount of load shedding and RE curtailment will not exceed $VaR_{\alpha}(\tilde{\mathfrak{n}}_{t})$ with probability $\alpha$,
\begin{align}\label{eq:eighteen}
VaR_{\alpha}(\tilde{\mathfrak{n}}_{t})=\min \{z| F_{\tilde{\mathfrak{n}}_{t}}(z) \geq \alpha \}.
\end{align}
$VaR_{\alpha}(\tilde{\mathfrak{n}}_{t})$ is a measure of risk, however it is not a reliable measure if $\tilde{\mathfrak{n}}_{t}$ has a fat tail distribution. $VaR_{\alpha}(\tilde{\mathfrak{n}}_{t})$ provides no information about the amount of $\tilde{\mathfrak{n}}_{t}$ that may occur beyond the value indicated by this measure \cite{Rockafellar}. CVaR is defined as the conditional expectation of load shedding and RE curtailment above the amount VaR$_\alpha$. Let $E$ denote the expectation over $n_t$.
\begin{align}\label{eq:ninteen}
CVaR_{\alpha}(\tilde{\mathfrak{n}}_{t})=E [ \tilde{\mathfrak{n}}_{t} | \tilde{\mathfrak{n}}_{t} > VaR_{\alpha}(\tilde{\mathfrak{n}}_{t}) ],
\end{align}
\begin{align}\label{eq:twenty}
CVaR_{\alpha}(\tilde{\mathfrak{n}}_{t})=\int_{-\infty}^{\infty}z dF_{\tilde{\mathfrak{n}}_{t}}^{\alpha}(z),
\end{align}
where
\[
    F_{\tilde{\mathfrak{n}}_{t}}^{\alpha}(z)=
\begin{cases}
    0,& \text{if }  z < VaR_{\alpha}(\tilde{\mathfrak{n}}_{t}) \\
    \frac{F_{\tilde{\mathfrak{n}_t}}(z)-\alpha}{1-\alpha}, & \text{otherwise}
\end{cases}.
\]
$CVaR_{\alpha}(\tilde{\mathfrak{n}}_{t})$ quantifies the value of the tail distribution of $\tilde{\mathfrak{n}}_{t}$ beyond the value of $VaR_{\alpha}(\tilde{\mathfrak{n}}_{t})$. $CVaR$ is a more conservative measure of risk than $VaR$, which is a lower bound on the risk. The function $g_{t}$ is defined as
\begin{align}\label{eq:objfunction}
g_{t}(b_{t})=   CVaR_{\alpha}(\tilde{\mathfrak{n}}_{t}),
\end{align}
for every $t \in \{1,...,T\}$. At each given time $t$, the energy storage unit minimizes its own objective function on $b_t,...,b_T$ as follows
\begin{align}\label{eq:obj100}
\underset{w.r.t. \,\ (\ref{eq:one}), (\ref{eq:constraint})}{\min\limits_{b_t,...,b_T}} \,\  \sum_{\tau=t}^{T} g_{\tau}(b_{\tau}).
\end{align}
Let $f_{t}$ be the density function of $n_t$. We assume the density function of ${n}_t$ satisfies Assumption~$1$.\\
\\
\textbf{Assumption~$1$:} $f_{t}$ is strictly positive on the interval $[P^{\min} \,\ P^{\max}]$.\\
\\
In the next section, it is proven that if Assumption~$1$ holds, then the optimal decision for the energy storage with objective function (\ref{eq:objfunction}) and (\ref{eq:obj100}) is unique for every $t \in \{1,...,T\}$.

\section{DYNAMIC PROGRAMMING FORMULATION}\label{SectionIII}

In this section, the optimization problem in (\ref{eq:obj100}) is reformulated as a dynamic programming problem. Let $J_{t}$ be the optimal objective function for the $(T-t)$-stage problem that starts at state $s_{t}$ at time $t$, and ends at time $T$,
\begin{align}\label{eq:obj11}
J_{t}(s_{t})=   \underset{w.r.t. \,\ (\ref{eq:one}), (\ref{eq:constraint})}{\min_{b_{t},...,b_{T}}} \sum_{\tau=t}^{T} g_{\tau}(b_{\tau}).
\end{align}
By applying the Dynamic Programming (DP),
\begin{align}\label{eq:DP}
& J_{T+1}(s_{T+1})=0, \notag \\
& J_{t}(s_{t})=  \underset{w.r.t. \,\ (\ref{eq:one}), (\ref{eq:constraint})}{\min_{b_t}} \{g_{t}(b_{t})+J_{t+1}(s_{t+1})\}.
\end{align}
By introducing the function $G_{t}$ for the given $s_t$
\begin{align}\label{eq:obj111}
G_{t}(b_{t}|s_t)=   g_{t}(b_{t})+J_{t+1}(a s_{t}+b_{t} \Delta t),
\end{align}
the DP equation (\ref{eq:DP}) can be written as
\begin{align}\label{eq:obj1111}
J_{t}(s_{t})=  \min_{   w.r.t \,\ (\ref{eq:one}), (\ref{eq:constraint})   }  G_{t}(b_{t}|s_t),
\end{align}
where
\begin{align}\label{eq:initialstate}
s_1:=\arg \min_{s^{\min} \leq s_{1} \leq s^{\max}} J_{1}(s_{1}).
\end{align}
\\
%In the following Proposition, it is shown that (\ref{eq:objfunction}) is strictly convex in $b_t$.
While the convexity of $CVaR_{\alpha}(\tilde{\mathfrak{n}}_{t})$ in $\tilde{\mathfrak{n}}_{t}$ is known \cite{Rockafellar}, in the following Proposition we prove the \textbf{strict} convexity of $CVaR_{\alpha}(\tilde{\mathfrak{n}}_{t}(b_t))$ in $b_t$ if Assumption~$1$ holds.\\
\\
\textbf{Proposition~$1$:} Let Assumption~$1$ hold. Then, $g_t(b_t)$ is strictly convex in $b_t$ for all $t=1,...,T$.

\begin{proof}The proof is given in Appendix A.\end{proof}

The constraints in minimization (\ref{eq:obj1111}) interconnect the minimization of the first and second term in (\ref{eq:obj111}). In the following Theorems, it is shown that (\ref{eq:obj1111}) is a convex optimization problem and therefore can be solved efficiently, e.g. by interior point methods \cite{Boyd}.\\
\\
\textbf{Theorem~$1$:} Let Assumption~$1$ hold. Then $G_{t}(b_{t}|s_{t})$ is convex in $(b_{t},s_{t})$ for the given $s^{\min} \leq s_{t} \leq s^{\max}$, and for all $b_t$ that satisfies (\ref{eq:constraint}).

\begin{proof} $J_{T+1}$ is the zero function and because of the convexity of $g_{T}(b_{T})$ in $b_{T}$ (Proposition~$1$), $G_{T}(b_{T}|s_T)$ is convex in $(b_{T},s_{T})$. Assume that $G_{t}(b_{t}|s_{t})$ is convex in $(b_{t},s_{t})$, the convexity of $G_{t-1}(b_{t-1}|s_{t-1})$ in $(b_{t-1},s_{t-1})$ is shown below. Let
\begin{align}\label{eq:obj11111}
G_{t-1}(b_{t-1}|s_{t-1})=   g_{t-1}(b_{t-1})+J_{t}(s_{t}).
\end{align}
In Theorem~$2$, it is proven that $J_t(s_t)$ is convex in $s_t$ if $G_t(b_t|s_t)$ is convex in $(b_t,s_t)$. Therefore from (\ref{eq:obj11111}) and Proposition~$1$, $G_{t-1}(b_{t-1}|s_{t-1})$ is convex in $(b_{t-1},s_{t-1})$.
\end{proof}

\textbf{Theorem~$2$:} Let $s^{\min} \leq s_{t} \leq s^{\max}$ and $b_t$ satisfies inequality (\ref{eq:constraint}). If $G_{t}(b_{t}|s_{t})$ is convex in $(b_{t},s_{t})$ then $J_{t}(s_{t})$ is convex in $s_{t}$.

\begin{proof} The proof is given in Appendix B. \end{proof}

The following Corollary implies the uniqueness of the optimal decision for the energy storage.\\
\\
\textbf{Corollary~$1$:} If Assumption~$1$ holds, then for the given $s^{\min} \leq s_{t} \leq s^{\max}$, $G_{t}(b_{t}|s_{t})$ is strictly convex in $b_t$ for all $b_t$ that satisfies (\ref{eq:constraint}).

\begin{proof} From Proposition~$1$, $g_t(b_t)$ is strictly convex in $b_t$ and from Theorem~$1$ and Theorem~$2$, $J_{t+1}(s_{t+1})$ is convex in $s_{t+1}$. By induction and (\ref{eq:obj111}), it can be concluded that $G_{t}(b_t|s_t)$ is strictly convex in $b_{t}$ for all $t=1,...,T$. \end{proof}

The proofs of Theorem ~$1$, Theorem~$2$, and Corollary~$1$ hold if we incorporate energy loss during the charging/discharging process of energy storage for a more realistic model. For example, the actual change of state of charge is $\eta_{in} b_t \Delta t$ when $b_t >0$, and $\frac{1}{\eta_{out}}  b_t \Delta t$ when $b_t <0$, where $0< \eta_{in}, \eta_{out}< 1$ are the efficiency factors. Several distributions (e.g. Weibull, normal, Erlang, and beta) have been used to model the variations in the electric load \cite{Singh}-\cite{Heunis}. Regardless of the type of distribution, the theoretical result of this work remains valid as long as the density function of $n_t$ satisfies Assumption~$1$. Additionally, the proofs of Theorem~$1$, Theorem~$2$, and Corollary~$1$ remain valid if the stage cost function, $g_t(b_t)$, is any strict convex function in the charging/discharging rate of the battery. We use a normal distribution model for the $n_t$ in the following section for illustration.
%The same analysis still hold, with this model.
\section{SIMULATIONS} \label{SectionV}
\begin{figure}[h!]
\includegraphics[width=9.5cm]{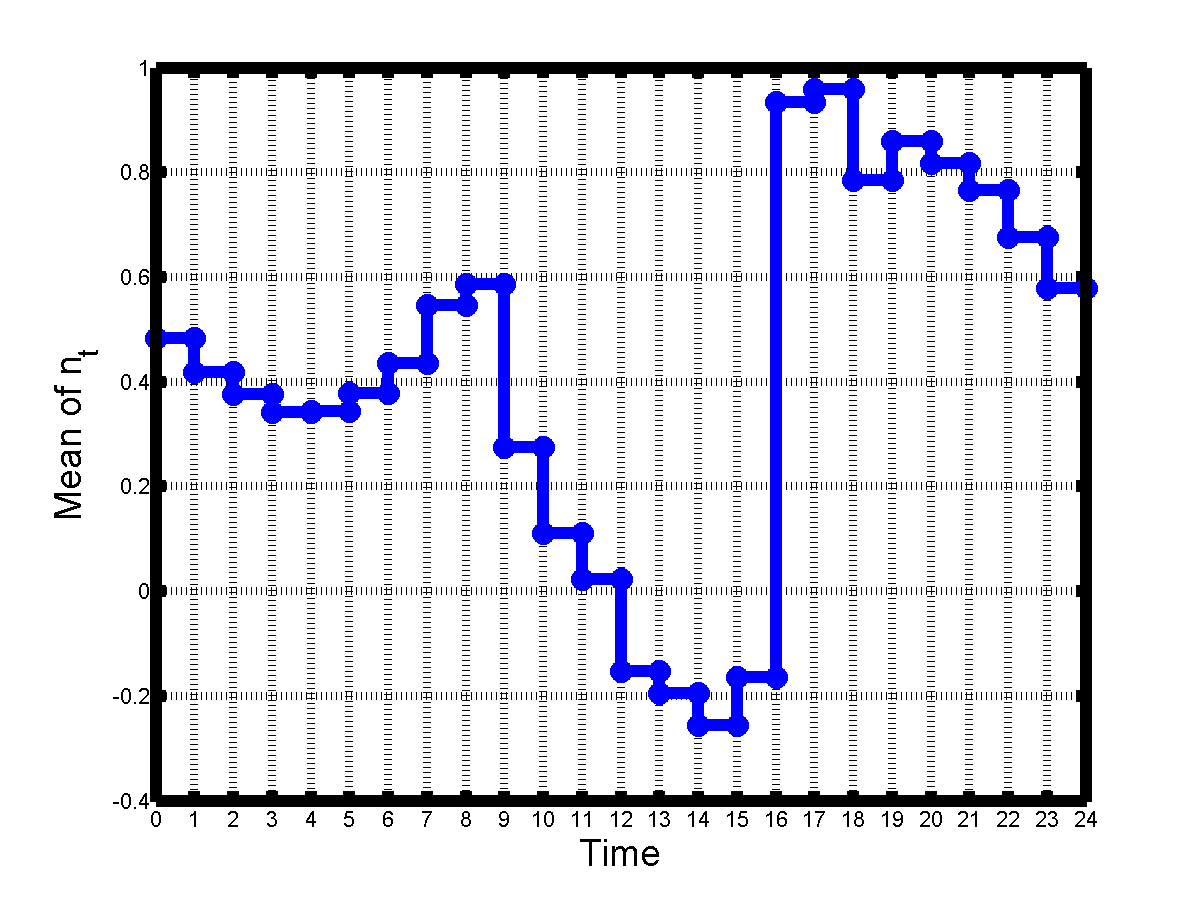}
\caption{Mean of $n_t$.}
\label{fig:NT}
\end{figure}

\textbf{Example:} We consider a $24$-stage storage optimization problem with given normalized parameters: $P^{\max}=0.6$, $P^{\min}=0$, $s^{\max}=1$, $s^{\min}=0$, $\alpha=0.01$ and $\Delta t=1$ hour. The random variable $n_t$ for all $t=0,...,23$ has a gaussian distribution. The mean of $n_t$ for all $t=0,...,23$ is plotted in the Figure~\ref{fig:NT}. The standard deviation of $n_t$ is $0.25$.

\begin{figure}[h!]
\includegraphics[width=9.5cm]{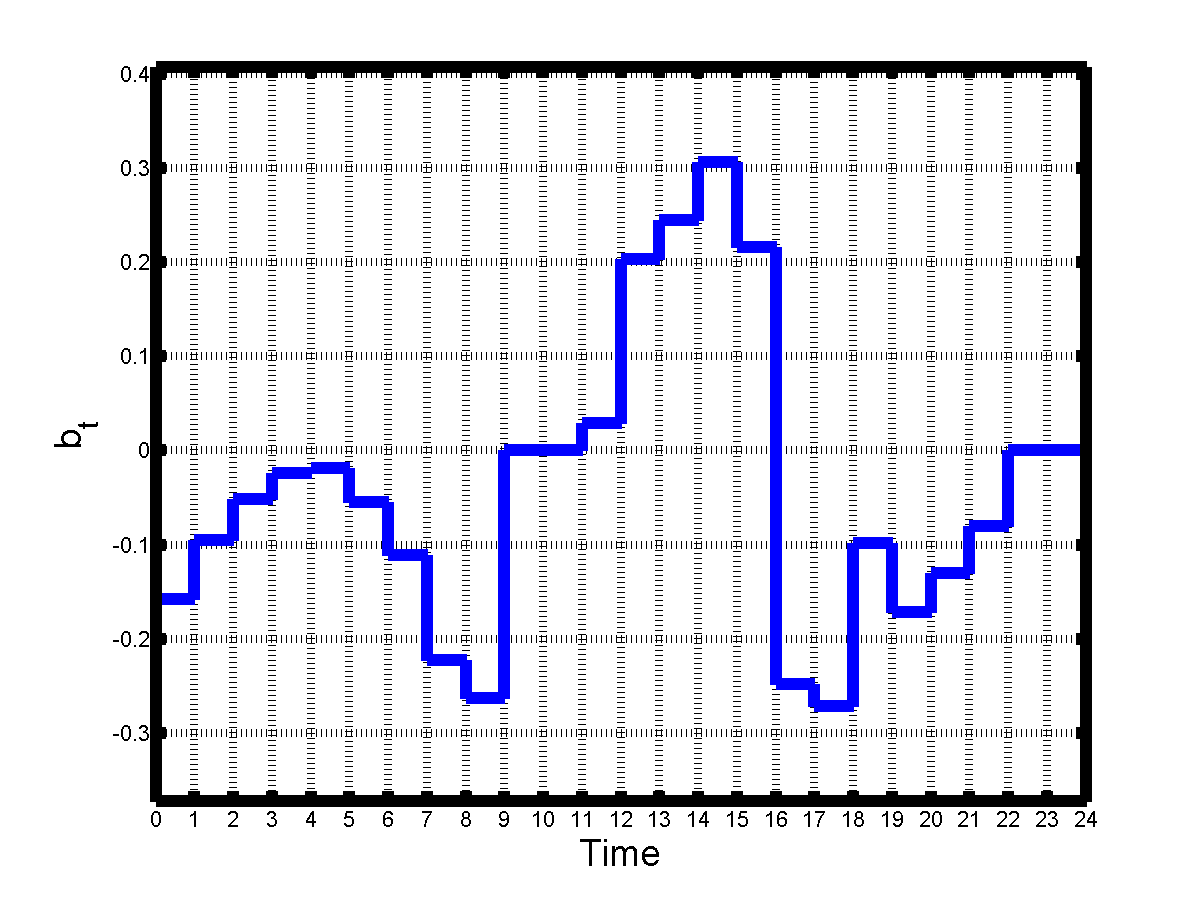}
\caption{The charge and discharge rate of energy storage.}
\label{fig:BT}
\end{figure}

\begin{figure}[h!]
\includegraphics[width=9.5cm]{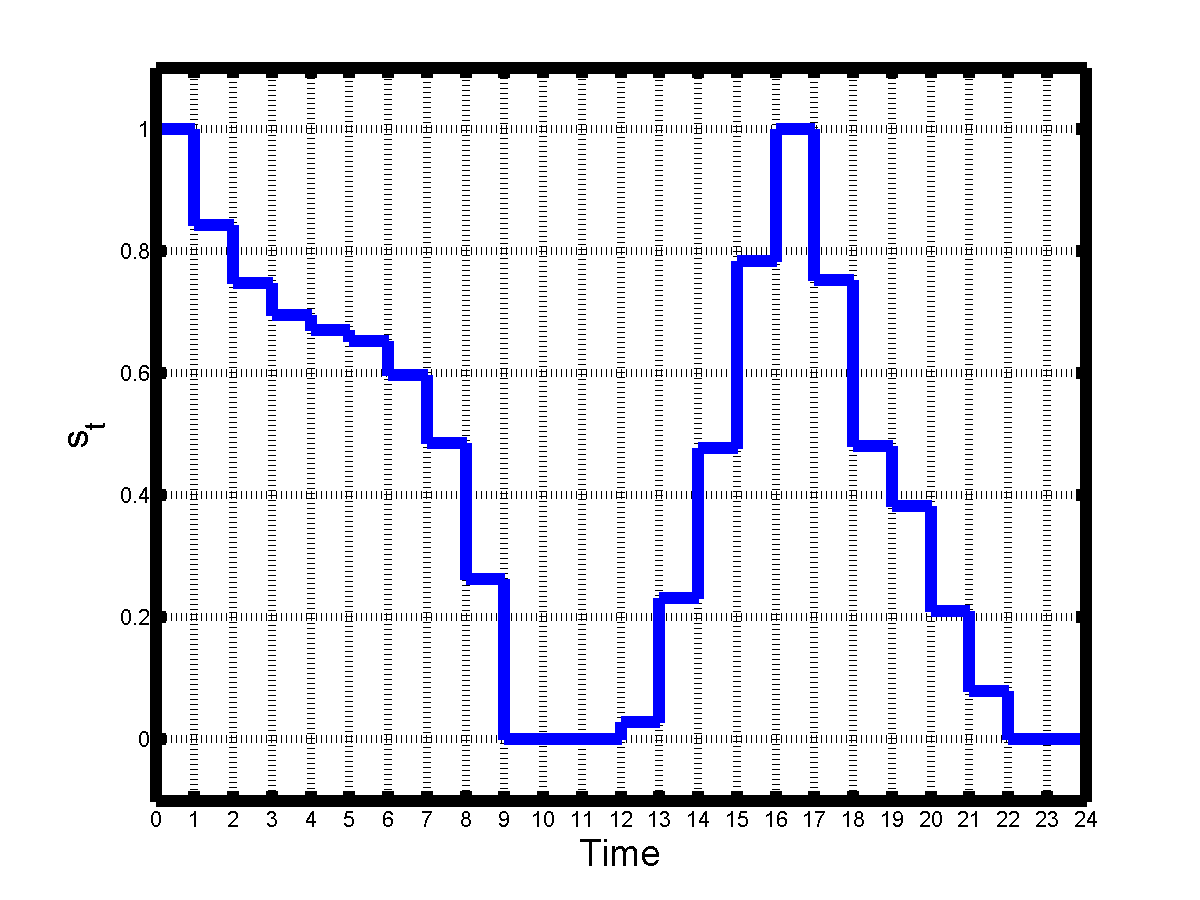}
\caption{State of charge of battery, $s_t$.}
\label{fig:SOC}
\end{figure}

For the sake of simplicity, the presented simulations are limited to $24$ hours. However, in practice the analysis should be over the lifetime of the energy storage. We consider a realization of $n_t$, that is equal to the mean of the random variable $n_t$ as plotted in Figure~\ref{fig:NT}. The charge and discharge rate of the energy storage is plotted in Figure~\ref{fig:BT}. The energy level of the battery is plotted in Figure~\ref{fig:SOC}.

\begin{figure}[h!]
\includegraphics[width=9.5cm]{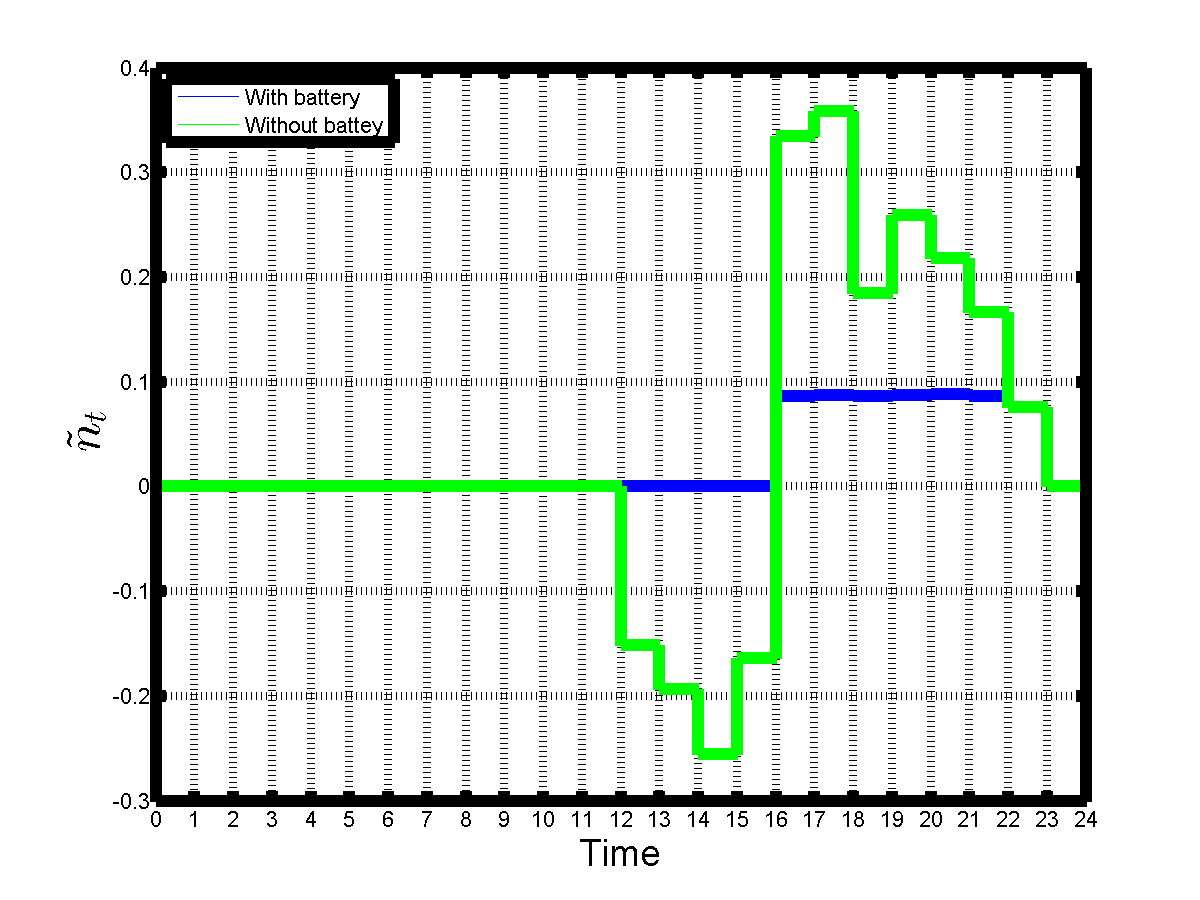}
\caption{Load shedding (positive value) and renewable energy curtailment (negative value) with battery (blue) and without battery (green).}
\label{fig:NTILDE}
\end{figure}

It is observed from Figures \ref{fig:BT} and \ref{fig:SOC} that the battery is charged during the daytime with excess renewable energy generation and discharged during the morning and evening peak times. The load shedding and renewable energy curtailment for both scenarios, with battery (blue) and without battery (green), are plotted in Figure~\ref{fig:NTILDE}. It is observed that the energy storage with an optimal strategy, decreases renewable energy curtailment and load shedding significantly.

\begin{figure}[h!]
\includegraphics[width=9.5cm]{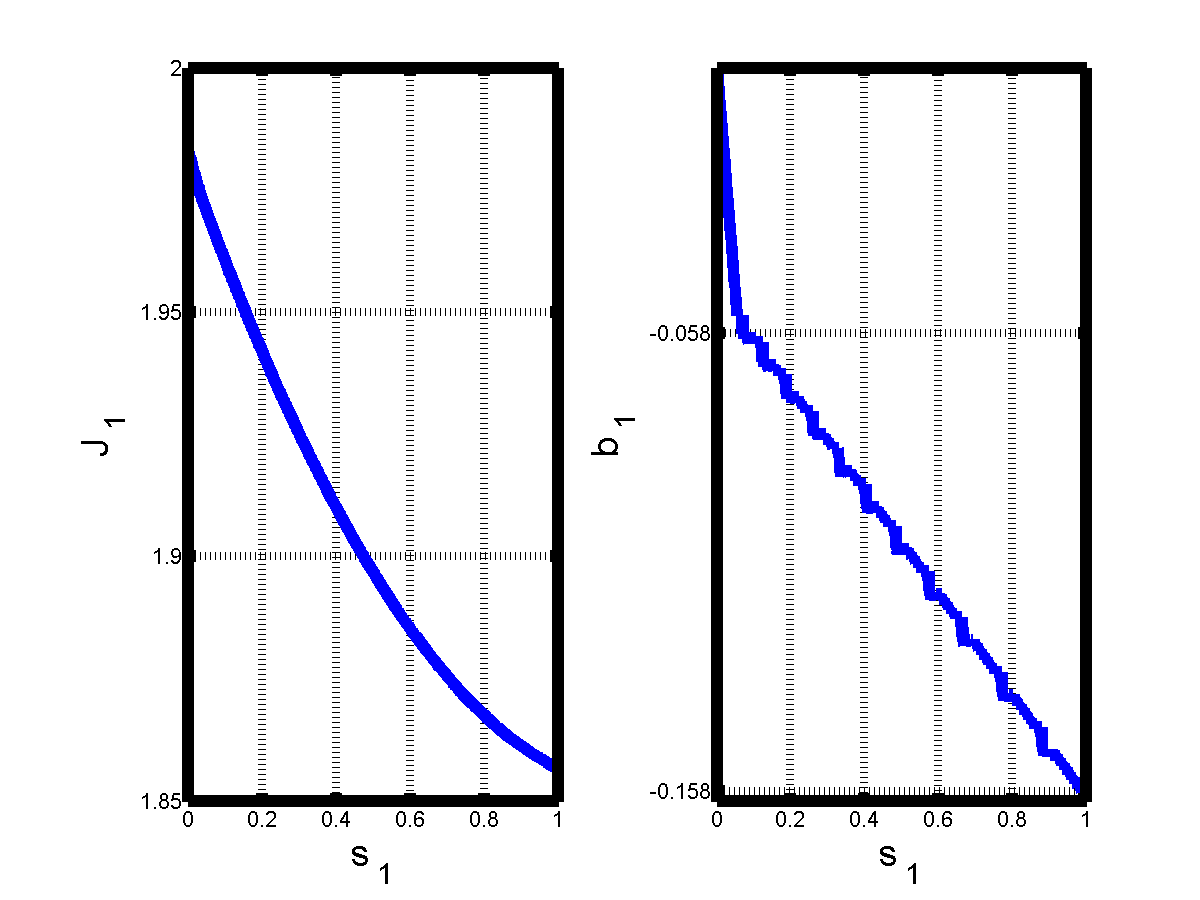}
\caption{The optimal objective function and charge/discharge rate from 12:00 A.M. to 1:00 A.M.}
\label{fig:stage1}
\end{figure}

\begin{figure}[h!]
\includegraphics[width=9.5cm]{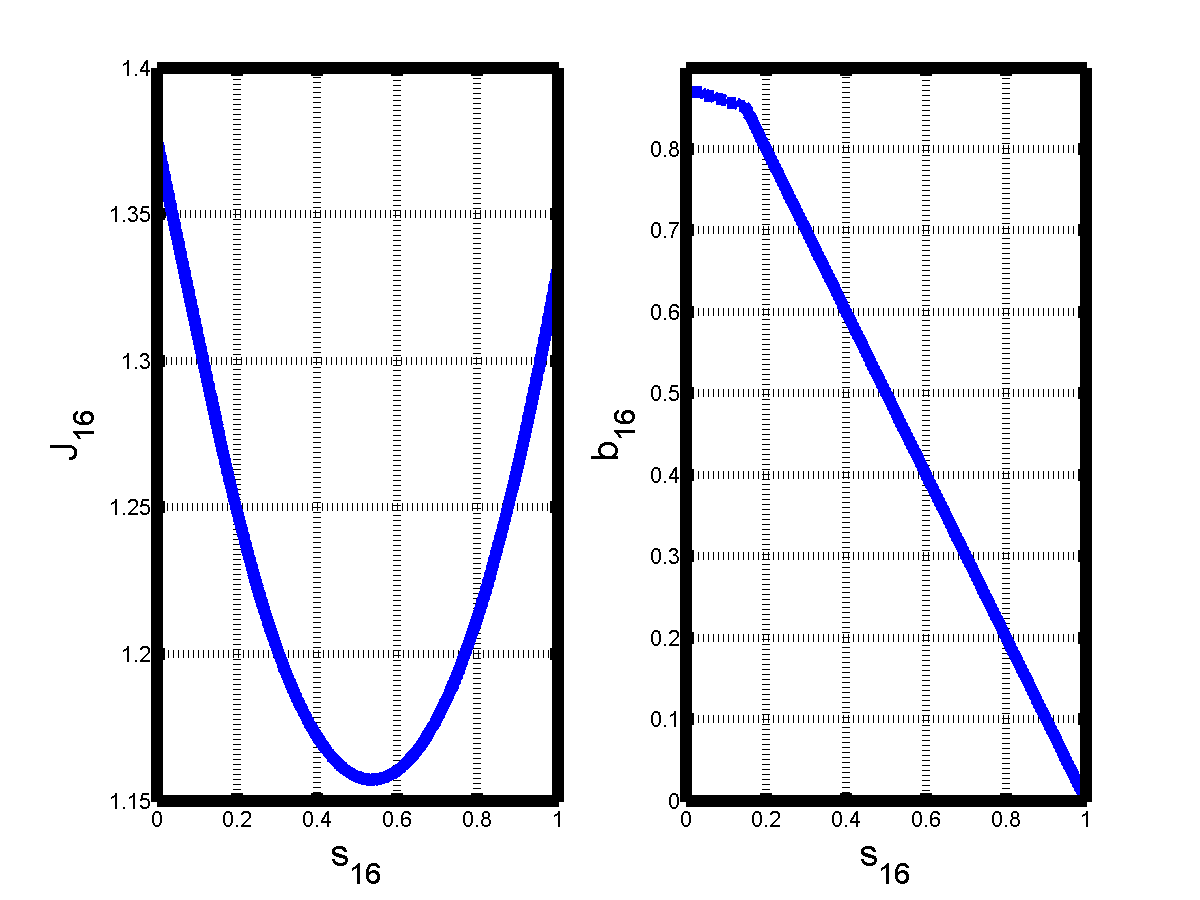}
\caption{The optimal objective function and charge/discharge rate from 3:00 P.M. to 4:00 P.M.}
\label{fig:stage16}
\end{figure}

In Figures~\ref{fig:stage1}-\ref{fig:stage16}, the optimal objective function and charge/discharge rate is plotted as a function of the energy level of the energy storage for the time slots $12:00$ A.M.$-1:00$ A.M. and $3:00$ P.M.$-4:00$ P.M. The battery is charged from $3:00$ P.M. to $4:00$ P.M. and discharged from $12:00$ A.M. to $1:00$ A.M. It is observed that the optimal objective function $J_t(s_t)$ is convex in the state of charge of the battery $s_t$. It is observed from Figure~\ref{fig:stage1} that the optimal energy level for the battery is $1.00$, which is assumed to be the initial energy level of the battery in the Figure~\ref{fig:SOC}.

\newpage
\section{CONCLUSIONS} \label{SectionVI}
In this work, a microgrid with renewable energy resources and an energy storage unit is considered. It is assumed that the load and RE are random processes over a finite time horizon. The main grid provides the microgrid operator with a total net load trajectory in order to satisfy the grid constraints. The net load trajectory is satisfied by the microgrid controller by applying renewable energy curtailment and load shedding. The energy storage unit is controlled by the microgrid independently from the main grid controller. The microgrid controller decides on charging and discharging strategies in order to minimize the risk of high amounts of load shedding and RE curtailment within a finite time horizon. The energy storage objective function in the grid is formulated as a finite stage dynamic programming problem. It is proven that under sufficient conditions the multi-stage objective function of energy storage is \textbf{strictly} convex in the charging/discharging rate. Thus, the optimal decision at each stage is unique for the given state of charge of the battery. The convexity of the objective function in charge/discharge rate, is needed for implementing mathematical programming that searches within the action space. Finally, numerical results are presented for an energy storage unit with a multistage objective function in a microgrid with random load and RE. It is observed that an optimal charging/discharching strategy for the energy storage based on an stochastic model for the load and renewable generation, can decrease the realized RE curtailment and load shedding.

%In order to mitigate the intermittencies of RE and relieve the grid of congestion and of short-supply and/during high-demand situations, it is envisioned that future smart grid systems will provide incentives for consumers to be more flexible in their consumption behavior. The potential for load control in electric grids is significant, and it is reasonable to assume that a portion of the load is flexible and can be deferred to the future load. There is a need to generate heuristic load shedding and renewable resource curtailment strategies based on the characteristics of stochastic loads, RE, and the dynamics of the grid. A poorly designed load-shedding strategy can worsen the system and cause/result in disruption.

\section*{APPENDIX}
\subsection{Proof of Proposition~$1$}
\begin{proof}
Let $f_t$ be the density function of $n_t$. Without loss of generality assume $\alpha=0$, from (\ref{eq:eighteen})-(\ref{eq:objfunction})
\begin{align}
g_t(b_t)= & \int_{P^{\max}-b_t}^{\infty} (n_t+b_t-P^{\max}) f_t \partial n_t\\ \notag
&- \int_{-\infty}^{P^{\min}-b_t}(n_t+b_t-P^{\min}) f_t \partial n_t.
\end{align}
Below, it is shown that $\int_{P^{\max}-b_t}^{\infty} (n_t+b_t-P^{\max}) f_t \partial n_t$ is strictly convex in $b_t$. Without loss of generality, assume $x < \gamma x+ (1-\gamma)y < y$. Below it is shown that for all $\gamma \in (0,1)$ and $x \neq y$
\begin{align}\label{ineq:con}
&\int_{P^{\max}-\gamma x-(1-\gamma)y}^{\infty} \big (n_t+\gamma x+(1-\gamma)y-P^{\max} \big ) f_t \partial n_t & \\ \notag
& \,\ \,\ \,\ \,\ \,\ \,\ \,\ \,\ \,\ \,\  \,\ \,\ \,\ \,\ \,\  < & \\ \notag
&\gamma \int_{P^{\max}- x}^{\infty} (n_t+ x-P^{\max}) f_t \partial n_t& \\ \notag
&+(1-\gamma) \int_{P^{\max}-y}^{\infty} (n_t+y-P^{\max}) f_t \partial n_t.&
\end{align}
The right hand side of inequality (\ref{ineq:con}) is equal to
\begin{align}
&\gamma \int_{P^{\max}-\gamma x-(1-\gamma)y}^{\infty} (n_t+x-P^{\max}) f_t \partial n_t \\ \notag
&- \gamma \int_{P^{\max}-\gamma x-(1-\gamma)y}^{P^{\max}-x} (n_t+x-P^{\max}) f_t \partial n_t \\ \notag
&+(1-\gamma) \int_{P^{\max}- y}^{P^{\max}- \gamma x-(1-\gamma)y} (n_t+ y-P^{\max}) f_t \partial n_t \\ \notag
&+(1-\gamma) \int_{P^{\max}- \gamma x-(1-\gamma)y}^{\infty} (n_t+ y-P^{\max}) f_t \partial n_t \\ \notag
&=\int_{P^{\max}-\gamma x-(1-\gamma)y}^{\infty} \big (n_t+\gamma x+(1-\gamma)y-P^{\max} \big ) f_t \partial n_t \\ \notag
&-\gamma \int_{P^{\max}- \gamma x-(1-\gamma)y}^{P^{\max}- x} (n_t+ x-P^{\max}) f_t \partial n_t \\ \notag
&+(1-\gamma) \int_{P^{\max}-y}^{P^{\max}-\gamma x-(1-\gamma)y} (n_t+y-P^{\max}) f_t \partial n_t.
\end{align}
The strict convexity follows from Assumption~$1$, and
$$\int_{P^{\max}- \gamma x-(1-\gamma)y}^{P^{\max}- x} (n_t+ x-P^{\max}) f_t \partial n_t <0$$
and
$$\int_{P^{\max}-y}^{P^{\max}-\gamma x-(1-\gamma)y} (n_t+y-P^{\max}) f_t \partial n_t >0.$$
Similarly, it can be shown that $-\int_{-\infty}^{P^{\min}-b_t}(n_t+b_t-P^{\min}) f_t \partial n_t$ is strictly convex in $b_t$. Therefore, $g_t(b_t)$ is strictly convex in $b_t$. %\qed
\end{proof}


\subsection{Proof of Theorem~$2$}
\begin{proof}
Without loss of generality, assume
\begin{align}\label{eq:const1}
s^{\min} \leq x < \alpha x+ (1-\alpha) \acute{x} < \acute{x} \leq s^{\max},
\end{align}
It is shown below that for all $\alpha \in (0,1)$, and $x$ and $\acute{x}$ that satisfy (\ref{eq:const1})
\begin{align}\label{eq:concavity1}
J_{t}\big(\alpha x+(1-\alpha) \acute{x} \big) \leq  \alpha J_{t}(x)&+(1-\alpha) J_{t}(\acute{x}).
\end{align}
Suppose, there exists an $x$ and $\acute{x}$, $s^{\min} \leq x < \acute{x} \leq s^{\max}$, and $\alpha \in (0,1)$ such that
\begin{align}\label{eq:concavity2}
J_{t}\big(\alpha x+(1-\alpha) \acute{x} \big)  >  \alpha J_{t}(x)+(1-\alpha) J_{t}(\acute{x}).
\end{align}
Let
\begin{align}\label{eq:concavity3}
b_{t}^m (x)= {\arg {\min_{ { \frac{s^{\text{min}}-a x}{\Delta t} \leq b_{t} \leq \frac{s^{\text{max}}-a x}{\Delta t}} }} G_{t}(b_{t}|x)},
\end{align}
$$\bar{x}=\alpha x+(1-\alpha)\acute{x}.$$
It can be proven by induction that $G_{t}(b_t|x)$ is continuous in $b_{t}$ for all $t=1,...,T$. Because of the compactness of the domain in (\ref{eq:concavity3}) and the continuity of $G_{t}(b_t|x)$ in $b_{t}$ the minimizer exists. From (\ref{eq:obj1111}), (\ref{eq:concavity2}), and (\ref{eq:concavity3}), it can be concluded that
\begin{align}\label{eq:concavity4}
G_{t}(b_{t}^m(\bar{x})| \bar{x})  > \alpha G_{t}(b_{t}^m(x)|x) +  (1-\alpha) G_{t}(b_{t}^m(\acute{x})|\acute{x}).
\end{align}
It is evident that
\begin{align}\label{eq:concavity5}
\frac{s^{\text{min}}-\bar{x}}{\Delta t}  \geq \alpha b_{t}^m(x)+(1-\alpha) b_{t}^m(\acute{x}) \geq  \frac{s^{\text{max}}-\bar{x}}{\Delta t}.
\end{align}
From (\ref{eq:concavity3}), (\ref{eq:concavity4}), and (\ref{eq:concavity5}), it can be concluded that
\begin{align}\label{eq:concavity6}
& G_{t}(\alpha b_{t}^m(x)+(1-\alpha)b_{t}^m(\acute{x})|\bar{x}) \\ \notag
& \geq G_{t}(b_{t}^m(\bar{x})|\bar{x}) \\ \notag
& > \alpha G_{t}(b_{t}^m(x)|x) +  (1-\alpha) G_{t}(b_{t}^m(\acute{x})|\acute{x}).
\end{align}
The first inequality is because of (\ref{eq:concavity5}), and $b_{t}^m(\bar{x})$ being the minimizer of $G_{t}(.| \bar{x})$. The second inequality is because of (\ref{eq:concavity4}). Inequality (\ref{eq:concavity6}) contradicts the convexity of $G_t(b_t|s_t)$ in $(b_t,s_t)$, therefore $J_t(s_t)$ is convex for $s^{\max} \leq s_t \leq s^{\min}$.
\end{proof}

\end{document}